**Model selection for risk analysis of wastewater networks**

Aaron Dunton, PhD Candidate, Civil and Environmental Engineering, University of Illinois at Urbana-Champaign

**Abstract**

In this paper, we test various models of wastewater infrastructure for risk analysis and compare their performance. While many representations are available, existing studies do not consider selection of the appropriate model for risk analysis. In this paper, we define two characteristics of wastewater models: the network granularity and the fidelity of the governing equations. We consider different combinations of these characteristics to determine 6 network representations that could be used as the foundation for risk analysis. We test the performance of each model as compared to predictions from the most detailed model, the full network with dynamic wave flow equations. We demonstrate the model selection for Seaside, Oregon. We conclude that the full network granularity is needed as compared to a coarse network representation. For the fidelity of the governing equations, connectivity analysis is reasonable if the primary goal is to determine the spatial distribution of hazard impact. To more accurately predict nodal performance measures, the dynamic wave equations are needed as they capture important physical phenomena.

**1. Introduction**

Critical infrastructure, including wastewater infrastructure, are vital to the functioning and prosperity of modern societies (Gardoni and Murphy 2020). However, past events have demonstrated that these infrastructure are vulnerable to damage from natural and anthropogenic hazards (Todd et al. 1994, Liu et al. 2013, Kazama and Noda 2012). Recent research has focused on how to improve the performance of critical infrastructure in the aftermath of a disrupting event. Mitigation measures and recovery optimization can decrease disaster impact. Because these events are rare and uncertain, there is a need for modeling to support mitigation planning and recovery optimization. This paper focuses on how wastewater infrastructure should be modeled for risk analysis. In Section 2, we review wastewater risk analysis research with a focus on how the infrastructure is modeled. In general, different representations of wastewater infrastructure are available, but none of the reviewed studies properly considered model selection. Especially in the practical context where the information needed to define the model may not be available and may need to be collected or generated (e.g., Dunton and Gardoni 2024, Dunton and Gardoni 2023), selection of an appropriate model among the available options allows for efficient risk analysis.

In this paper, we test the performance of various models of wastewater infrastructure for risk analysis. We focus on networks, a specific type of model that describes the connection of components and subsequent flow dynamics. We start by defining network characteristics that have primary influence over the model running time. Different combinations of these characteristics define the 6 models that we consider. Following the steps for risk analysis, we then simulate damage to each of the components of the infrastructure. For this paper, we assume probabilities for each of the damage states, so general conclusions about network selection are applicable to any hazard which could damage the components of the wastewater network. Next, we assess the functionality of the infrastructure given the damage. For each model, we define a performance measure that describes the functionality of each node in the network. We compare the performance of the models by comparing the predictions of the node-level performance measure from each model with the predictions from the most detailed model. We demonstrate the model selection for risk analysis of the wastewater network of Seaside, Oregon. We estimate the expected value of the performance measure for 5,545 nodes for each of the 6 models given the probability of damage for each component. Finally, we compare the predictions from the 6 models.

The rest of the manuscript is organized as follows. In Section 2, we review previous research in risk analysis for wastewater infrastructure. In Section 3, we describe the network representations of wastewater infrastructure that we consider. In Section 4, we describe how we incorporate damage to the components. In Section 5, we describe the node-level performance measures that we use to assess the functionality of the network given the damage. In Section 6, we present the case study results. Finally, we conclude with Section 7.

## 2. Literature review

Here, we review the literature on risk analysis for wastewater networks, with particular emphasis on the network representation used in these studies. Some of the studies in risk analysis for wastewater systems have used system models that include subcatchments, tanks, pumps, combined sewer overflows, and treatment plant processes (Saagi et al. 2016, Sweetapple et al. 2018, Sweetapple et al. 2022). These models have been schematic representations of the system, lacking geographical information. The models are at a coarse granularity; sewer inflow is represented at the subcatchment scale, and the model represents how the subcatchments are connected with tanks, pumps, and the treatment plant. The treatment plant processes are modeled in detail. They study the effect of failure of different components of the system, including increased inflow at the subcatchment scale, failure of pipes, and failure of specific treatment processes.

Others have used network representations based on the geographic layout of wastewater networks. Rodriguez et al. (2013) and Sadashiva et al. (2021) developed network models for specific case study locations. Rodriguez et al. (2013) used a coarse representation of the network of Bogota, Colombia, where wastewater and stormwater inflow are at the subcatchment scale. Sewer trunks connect the subcatchments with the treatment plant, and combined sewer overflow locations are also included. For simulations, they proposed and used a conceptual model to predict pollutant properties based on monitoring data. Sadashiva et al. (2021) used a sewer network model provided from a local utility for four cities near Wellington, New Zealand. It is not clear whether the network was simplified, or if the full network information was available. For functionality assessment, they used connectivity analysis. They also considered whether the treatment plant was functioning, producing maps that showed the extent of three functionality states (no service; collected; collected and treated) over time. Other studies have used flow analysis (Zhang et al. 2017, Kamali et al. 2022) for functionality assessment, most commonly using the EPA's Storm Water Management Model (SWMM, EPA 2022).

Other studies focus on risk analysis for the treatment plant specifically. For modeling risk to sea-level rise and coastal flooding, the location of treatment plants and their elevation have been used (Hummel et al. 2018, Kool et al. 2020, Olyaei et al. 2018). More generally, Holloway et al. (2021) modeled the treatment plant processes and assessed the stress on the system from climate change, changes in wastewater inflow, and under-investment in the wastewater system.

Finally, many studies from asset management of wastewater networks have relied heavily on databases of component characteristics. Ugarelli et al. (2010a) emphasized the importance of these databases as the backbone for asset management. Many studies have focused on condition prediction and condition monitoring for damage to sewer components (Mohammadi et al. 2020, Bailey et al. 2015, Ugarelli et al. 2010b, Zamanian et al. 2021, Abebe and Tesfamariam 2020a, Davis et al. 2013). Others have studied how to handle incomplete information in these databases using data imputation (Kabir et al. 2020) and natural language processing of available records (Chahinian et al. 2021). For application to risk analysis, some have used these characteristic data to identify pipes with severe failure consequences (Abebe and Tesfamariam 2020b, Laakso et al. 2018, Balekelayi and Tesfamariam 2021). Some of these studies have used network characteristics (Laakso et al. 2018, Balekelayi and Tesfamariam 2021). That is, the underlying model used is a network that includes component characteristics.

In general, there are many models that have been used to represent wastewater infrastructure for risk analysis. Some models are from the more general asset management and operations research, while others are commonly used for hazard risk analysis. We emphasize the differences between the various reviewed models. First, different components of the wastewater system may or may not be modeled. For example, the treatment plant model may include individual treatment processes, it may be modeled using a simplified representation, or the treatment plant may not be modeled at all. Also, some studies used coarse representations of the infrastructure to represent the system-scale performance, while others used detailed pipe-level models. Second, different model fidelities are available. For example, some studies used flow analysis while other studies used topological metrics only. In general, previous studies only consider a single representation of wastewater infrastructure for risk analysis. They do not consider that many representations are available, and there may be some representation that is best for certain risk analyses. Addressing this limitation of existing research is the primary contributions of this paper.

## 3. Representations of wastewater infrastructure

Wastewater infrastructure provide sanitation by collecting and treating wastewater. In this paper, we focus on the collection system which transports wastewater from buildings to treatment facilities. Typical components of this system are shown in Figure 1. At the treatment plant, there would be an overflow from which untreated wastewater would be released if the treatment plant were not functional. So, the collection system provides the sanitation function in terms of preventing health hazards near occupied buildings. While there may be environmental impacts of combined sewer overflow that could be considered, we focus first on the main functioning of the wastewater system to move wastewater away from communities. Moreover, some amount of combined sewer overflow would be acceptable for the type of rare event that we are considering in this paper.

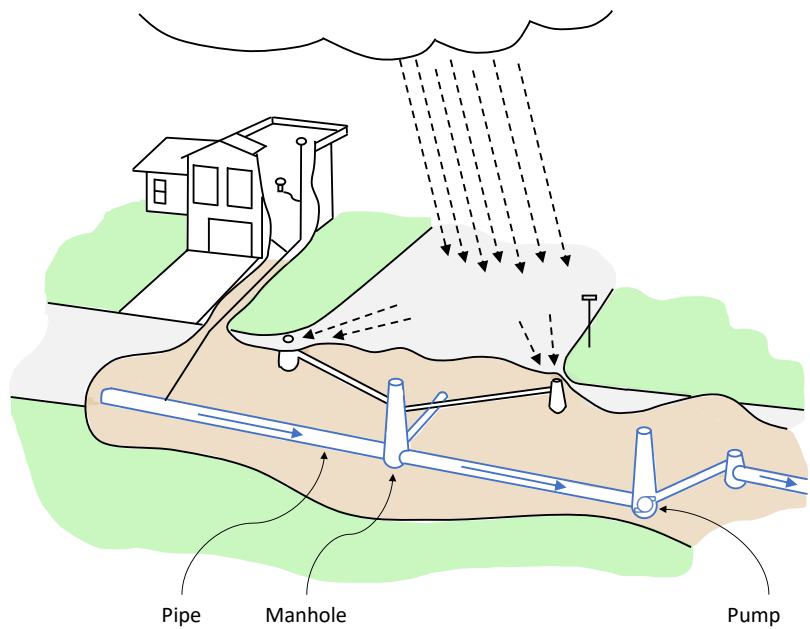

**Figure 1.** Wastewater collection system components

Networks are used to model wastewater collection systems. In this paper, we consider 6 unique combinations of two network characteristics to define a set of models. Each model represents the infrastructure at a different level of accuracy. In this paper, we are testing whether simpler models are accurate enough for risk analysis by comparing the results between simpler and more accurate models. The two characteristics that we consider are the primary contributions towards the model running time. That is, we consider a balance between accuracy and model running time, as is needed, e.g., for uncertainty quantification or optimization (i.e., applications where many model runs are needed). Table 1 lists the models that we test in this paper. The first characteristic that varies over the models is the network granularity. Network granularity is related to which components are included in the network. Previously, we defined a hierarchical conceptual model of infrastructure, visualized in Figure 2 (Dunton and Gardoni 2024). The arterial network includes the larger components that serve large areas of the community, and the capillary networks include the smaller components that connect to every building. For model selection, we consider two granularities: the arterial network and the full network (i.e., including both the arterial and capillary networks). For wastewater networks, the arterial network can be based on the location of pumping stations and clearly delineated basins from the tree topology.

**Table 1.** Combinations of network granularity and fidelity of the flow equations

|  | Dynamic Wave Equations | Kinematic Wave Equations | Connectivity |
|---|---|---|---|
| Full Network | Model DF | Model KF | Model CF |
| Arterial Network | Model DA | Model KA | Model CA |

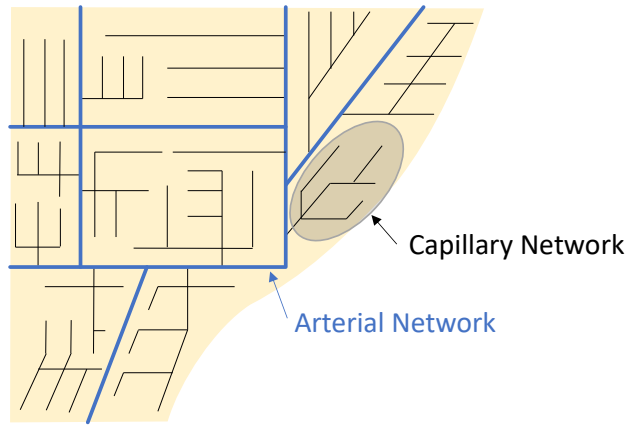

**Figure 2.** Hierarchical conceptual model of infrastructure, including the arterial and capillary networks

The second characteristic that varies between the models listed in Table 1 is the fidelity of the governing equations. The governing equations describe the network dynamics. At the lowest fidelity, we include connectivity analysis as an option for the governing equations. In connectivity analysis, undamaged components are represented as nodes and edges of a graph. If there is a path of undamaged components from a node to the treatment plant, then the node is connected. The connection or disconnection of each node is represented by an indicator variable with a value of 1 if the node is connected or a value of 0 if the node is disconnected.

While connectivity analysis only describes the possibility of flow from each node to the treatment plant, flow analysis additionally describes the actual movement of the wastewater. In flow analysis, wastewater inflow is produced at nodes. The amount of wastewater is commonly modeled with a baseline value and a time pattern, or, if available, measured wastewater flow data is used. Once wastewater is generated at nodes, it then flows along gravity-driven pipes. The gravity-driven flow of wastewater through a pipe is represented by the following dynamic wave equation:

$$\frac{\partial Q}{\partial t} + \frac{\partial (Q^2/A)}{\partial x} + gA\frac{\partial Y}{\partial x} = gA(S_0 - S_f) \tag{1}$$

where $Q$ is the flow, $t$ is time, $A$ is the flow area, $x$ is the location along the pipe, $Y$ is the flow depth, $S_0$ is the conduit slope, and $S_f$ is the friction slope (EPA 2017). After gravity-driven flow, wastewater may be stored in tanks that feed into pumps that are activated when the level of water in the tank reaches a certain value. The flow dynamics for each pump are represented by a pump curve. Next, the treatment plant processes may be modeled, and the flow is finally directed to the outfall. In this paper, we focus on modeling the collection system, so we do not model the treatment plant processes. Instead, flow is directly delivered to the outfall.

In the flow analysis, solving the dynamic wave equations for each pipe at each timestep is computationally demanding. The kinematic wave equation is the following simplified version of the dynamic wave equation:

$$0 = gA(S_0 - S_f) \tag{2}$$

The kinematic wave equation can be solved much faster than the dynamic wave equation, but it is less accurate. Specifically, the kinematic wave equation does not consider backwater or surcharge effects. To summarize, in this paper, we compare risk analysis results when using the following governing equations for the network dynamics, listed from highest to lowest fidelity: the dynamic wave equations, the kinematic wave equations, and connectivity analysis. These 3 options for fidelity are combined with the 2 options for granularity to define the 6 models that we test in this paper, listed in Table 1. However, many more combinations of granularity and fidelity can be defined by considering hybrid models that have multiple granularities or fidelities, over geographic extents or over the duration of a simulation.

**4. Damage prediction**

Numerous hazards may damage the components of the collection system. For example, past earthquakes have caused damage (Todd et al. 1994, Kazama and Noda 2012, Liu et al. 2013). Probabilistic damage models are used to predict the damage to components in future events as a function of the local hazard intensity measure, yielding the probability of each component being in each damage state. These damage models have conventionally been empirical models. For example, the models from the American Lifelines Alliance have commonly been used (ALA 2001). However, these types of models are based on limited historical data and do not incorporate the physics of the damage action. More detailed physics-based damage models can be used to evaluate the probability of damage to each component, incorporating component-specific characteristics (e.g., Iannacone and Gardoni 2023). In this paper, we do not focus on the damage models, but rather, selection of the representation of the network itself. To that end, we assume the probability that each component is in a certain damage state. For the purpose of selecting a network representation, assuming these probabilities is sufficiently realistic as both the conventional empirical damage models and the more detailed physics-based damage models produce these probabilities for each component. Moreover, by not considering a specific hazard, the general results from this analysis apply to any hazard that could damage the components. The damage probabilities that we use in this paper are listed in Table 2. The tanks and pumps are either undamaged or failed, each with a probability of failure of 0.1.

For the pipes, we define partial damage states based on a loss of flow capacity in addition to failure. We use simulations to determine how the uncertainty in the damage state of each component affects the uncertainty in the performance measure, the formulation of which we describe next.

Table 2. Assumed probability of damage for components

| Component | Damage State | Probability |
|---|---|---|
| Tanks | Undamaged | 0.9 |
| | Failed | 0.1 |
| Pumps | Undamaged | 0.9 |
| | Failed | 0.1 |
| Pipes | Undamaged | 0.85 |
| | 50% Flow Capacity | 0.1 |
| | 10% Flow Capacity | 0.04 |
| | Failed | 0.01 |

## 5. Functionality assessment

Finally, we assess the functionality of the system given the damage. We define a performance measure, $PM$, ranging from 0 to 1, and we evaluate its expected value for each demand node in the system given the damage probabilities from Section 4. If the value of the performance measure is 1 for a given node, then that node has not experienced any loss of functionality. This performance measure is the output of risk analysis. We consider model selection for risk analysis by comparing the performance measures for each demand node using the different representations.

For connectivity, the performance measure is the probability of disconnection of each node. That is, in individual simulations of the damage, each node is either disconnected ($PM = 0$) or connected ($PM = 1$). The probability of disconnection is then evaluated over the simulations. For flow analysis, the performance measure is again 0 when a node is disconnected from the treatment plant due to failure of pumps, tanks, or pipes. However, we additionally consider the influence of partial damage to the pipes. For nodes that are not disconnected, we define a performance measure, $PM_f$, that represents the ease with which water flows from each node to the treatment plant. For node $i$, the path to the treatment plant is the set of pipes $s_i$. For each simulation, the performance measure of node $i$, $PM_f^i$, is related to the number of pipes on the path with maximum flow greater than 80% of the flow capacity, as follows:

$$PM_f^i = 1 - \frac{\sum_{(j,k) \in s_i} I_{80}(j,k)}{n_i} \quad (3)$$

where $n_i$ is the number of pipes in $s_i$ and $I_{80}(j,k)$ is an indicator with value 1 if the maximum flow along pipe $(j,k)$ is greater than 80% of the pipe's capacity, and 0 otherwise. So, if only a single pipe in the path to the treatment plant is critically loaded, the value of the performance measure is nearly 1. This is a fair measure of performance for this type of system because if only one pipe in the path loses flow capacity, water would back up into preceding pipes and may cause pressurized flow in the critical pipe. The node itself might not see any difference in performance given these downstream effects. So, partial functionality in one pipe along the path is less critical than partial functionality in multiple pipes along the path, and $PM_f$ represents this. Moreover, this measure gives an indication of the amount of repair work that would be needed to restore functionality of each node. Again, in individual simulations, the performance measure is evaluated. Then, we take the sample average over the simulations as an estimate of the expected value of the performance measure for each node. The necessary number of simulations are determined according to a convergence criterion on the standard deviation of the estimate of the performance measure, $\delta_{PM}$, evaluated for each node as follows:

$$\delta_{PM} = \sqrt{\frac{1-PM}{N*PM}} \quad (4)$$

where $N$ is the number of simulations.

## 6. Case study

We demonstrate the procedure to test different models for the wastewater network of Seaside, Oregon, a testbed for community risk and resilience research (e.g., Guidotti et al. 2019, Boakye et al. 2019, Nocera and Gardoni 2022). Seaside is a small coastal community with a population that varies from approximately 7,000 to more than 20,000 during the summer tourist season. The location of Seaside and the city extents are show in Figure 3. The information needed to define the wastewater network are available, and the full and arterial networks are shown in Figure 4. To determine the arterial network, we include all tanks and pumps, and all nodes and pipes that are downstream from any tank or pump. To maintain the amount of wastewater inflow, the nodes immediately before the tanks/pumps are assigned all the wastewater inflow from the corresponding upstream nodes in the full network.

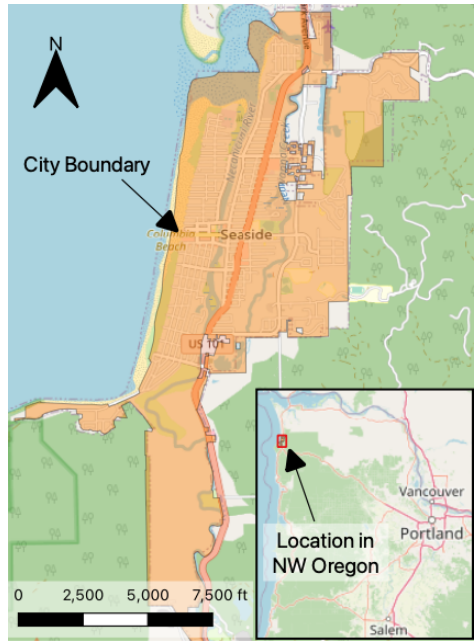

**Figure 3.** Case study location

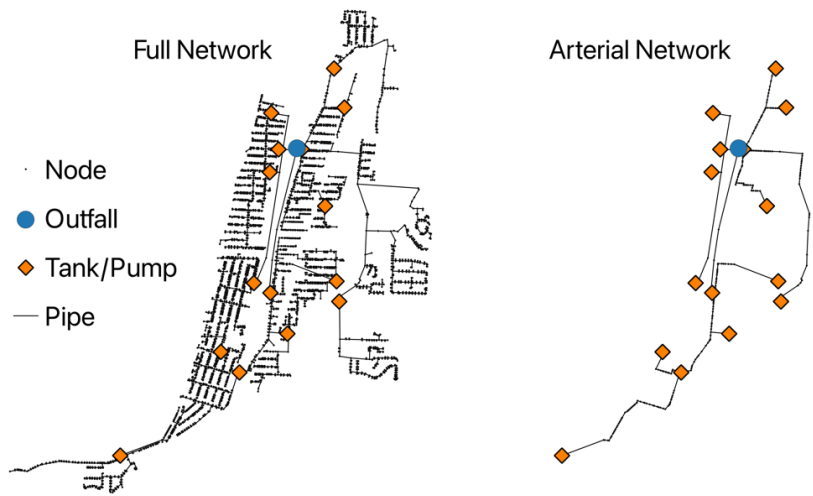

**Figure 4.** Full network and arterial network for Seaside

Next, we simulate the damage to the pumps, tanks, and pipes with the probabilities listed in Table 2. To evaluate the performance measure for each simulation, we use NetworkX for connectivity analysis (Hagberg et al. 2008), and we use the EPA's Storm Water Management Model (SWMM) for flow analysis (EPA 2022). In SWMM, we choose the corresponding routing equations (i.e., dynamic wave equations or kinematic wave equations) for each of the flow-based models. We also use GeoPandas to organize and save the node-level results as geospatial vector data (Jordahl et al. 2021). To determine the number of simulations, we use a convergence criteria. We evaluate $\delta_{PM}$ for each node after every 100 simulations. We continue until the average $\delta_{PM}$ over all of the nodes is less than 0.1. The convergence statistics for each of the models are reported in Table 3, along with the model running time. Overall, connectivity analysis took less than 10% of the time of flow analysis for both the full and arterial networks. For the full network, a simulation with kinematic wave routing took about 25% of the duration of a simulation with dynamic wave routing, as expected. However, for the smaller arterial network, a simulation with dynamic wave routing was slightly quicker than a simulation with kinematic wave routing.

**Table 3.** Model running times and convergence statistics

| Model Name | Number of Simulations | Average COV of PM | Time (s) | Time per Simulation (s) |
|---|---|---|---|---|
| DF | 300 | 0.0896 | 57935 | 193 |
| KF | 200 | 0.0690 | 9342 | 46.7 |
| CF | 100 | 0.0967 | 116 | 1.16 |
| DA | 100 | 0.0780 | 76.5 | 0.765 |
| KA | 100 | 0.0820 | 80.6 | 0.806 |
| CA | 100 | 0.0729 | 7.11 | 0.0711 |

The estimate of the expected value of the performance measure for every node is mapped in Figure 5 for the 6 models. For the arterial networks, the value of the performance measure for the capillary nodes that were not included in the arterial network is the value for the node in the arterial network to which the wastewater flow was directed. For models DA, KA, and CA, this yields areas with constant value of the performance measure. These areas come from the structure of the wastewater network and can also be seen in the results for models DF, KF, and CF. However, the performance measures from these models also have some local variations that the results from the arterial networks do not have. As compared to model DF, models KF and CF have similar local variations. Overall, the spatial distribution of impact on functionality is consistent between all of the models. That is, the areas with the worst and best performance are generally the same for every model. In general, the probability of damage increases as the distance from the outfall increases. That is, there are more components that may be damaged, so the probability of functionality loss increases. To summarize, the spatial distribution of the impact is consistently represented by all of the models, and all of the full networks capture similar local variations in the performance measure that are not captured by the arterial networks.

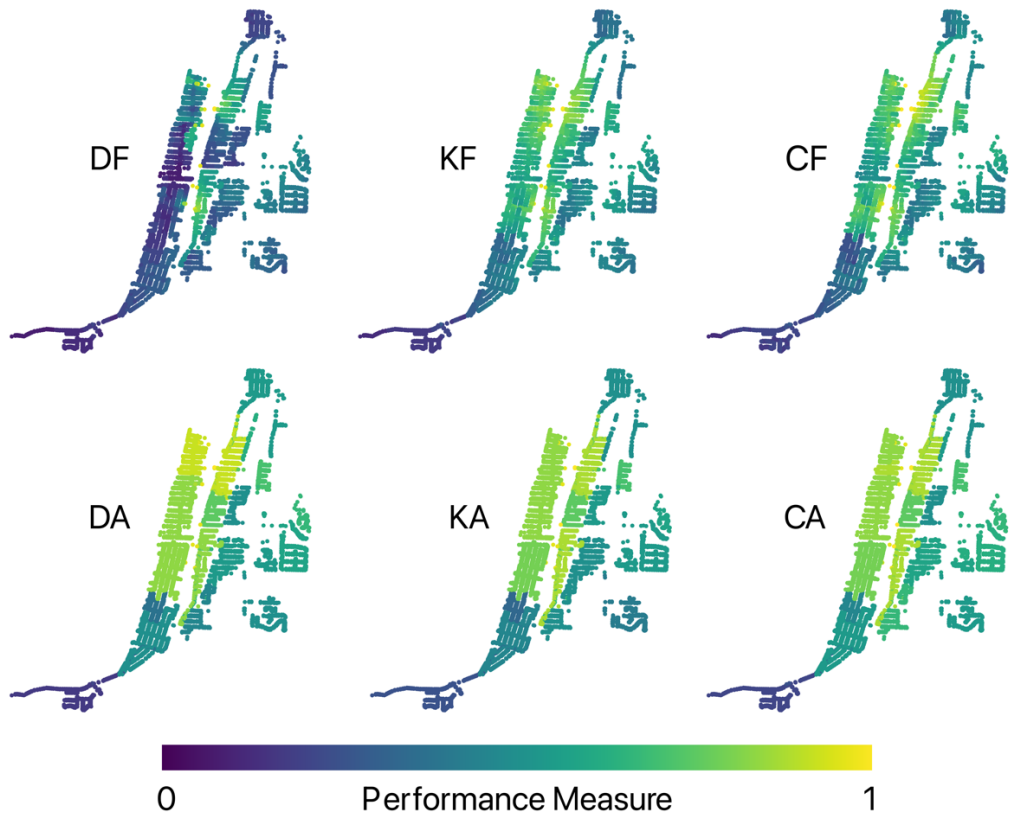

**Figure 5.** Maps of the performance measure for each model

While the spatial distribution of the performance measure is similar for all models, there are some significant differences in terms of value. To analyze these differences, Table 4 lists the average performance measure over the nodes for each model. Also, Figure 6 plots the performance measure for every node from each model versus the performance measure for every node from model DF (i.e., the most detailed model). We include the plot for model DF versus model DF to be consistent with Figure 5 and also to demonstrate that a perfect model should have points that follow the 1:1 line, also plotted. In general, model DF gives the lowest values of the performance measure. As compared to model KF, the flow analysis in model DF accounts for backwater and surcharge effects. That is, if a pipe loses flow capacity from damage and surcharges in the simulation, model DF is able to account for the effect on the flow in the preceding pipes where model KF is not. So, when partial damage to a pipe causes a loss of flow capacity, model DF will yield more pipes with flow greater than 80% of their flow capacity as compared to model KF. Since this is what $PM_f$ is counting, the values are lower for model DF. Similarly, model CF does not account for partial damage at all. In general, with increasing fidelity of the governing equations, partial damage is more accurately represented, resulting in lower performance measures. To compare the accuracy of the results using a single value, Figure 6 also includes the mean absolute error (MAE) for each model. Notably, model CF performs better than model DA. This indicates that, given a choice between finer network granularity or more accurate fidelity of the governing equations, finer network granularity would yield more accurate results for wastewater network risk analysis.

**Table 4.** Average nodal performance measure for each model

| Model | Average Performance Measure |
|:---:|:---:|
| DF | 0.336 |
| KF | 0.541 |
| CF | 0.546 |
| DA | 0.651 |
| KA | 0.627 |
| CA | 0.666 |

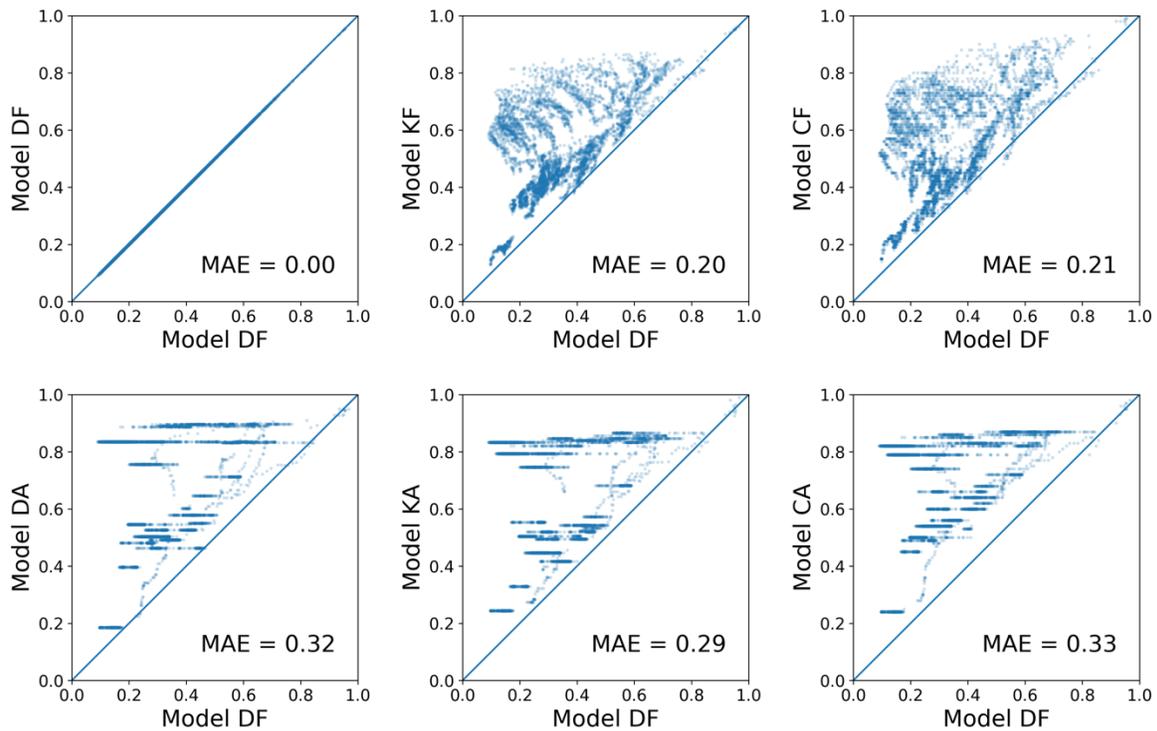

**Figure 6.** Nodal performance measures from each model versus performance measures from model DF

Because the results in Figure 6 seem to indicate similar performance between models KF and CF, and models KA and CA, we compare these models in Figure 7. Models CF and KF give similar predictions with a MAE of 0.04. The scatter in the relation is mostly due to the estimation error from the simulations. The MAE between model KA and CA is also low (0.05). These results indicate that flow analysis using the kinematic wave equations is not providing information about the functionality of wastewater networks that is not captured by connectivity analysis. To capture the relevant physics of the flow for risk analysis, the dynamic wave equations are needed.

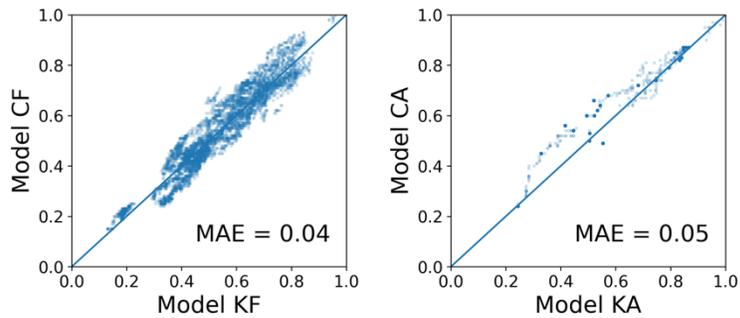

**Figure 7.** Comparison of nodal performance measures between models using kinematic wave routing and connectivity analysis

To directly assess the effect of using the full or arterial network, Figure 8 compares the results from models DF and DA; KF and KA; and CF and CA. For the kinematic wave routing and connectivity analysis, the full network gives slightly lower values of the performance measure than the arterial network. With more components to be damaged, there is more damage on average, yielding lower values of the performance measure for the full networks. With dynamic wave routing, the difference is more significant. This indicates that surcharging in the capillary pipes has significant effect on the value of the performance measure. To fully capture the relevant physics, the full network with dynamic wave routing is needed. That is, even with the highest fidelity of the governing equations, the arterial network does not accurately capture the salient physics. This is consistent with the previous observation that model CF performs better than model DA.

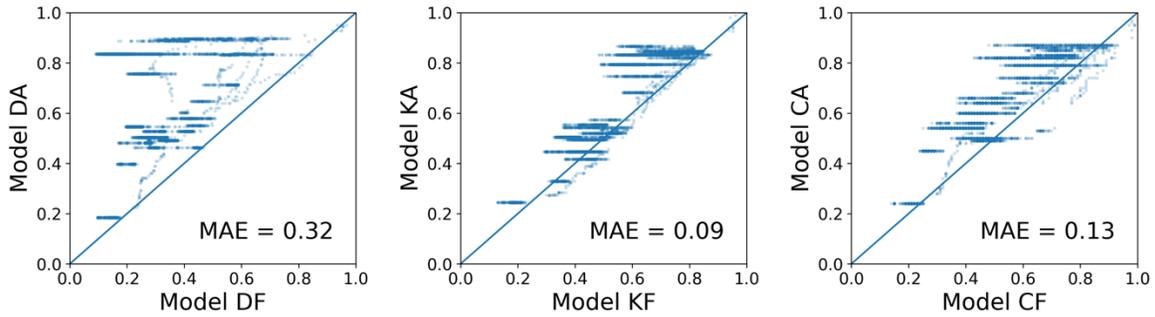

**Figure 8.** Comparison of nodal performance measures between models using the arterial network and the full network

Figure 9 shows a final comparison between the tested models in terms of their accuracy and the running time per simulation. First, all of the arterial models performed poorly, even when compared with model CF. While this could indicate that the full model should be used, we also acknowledge that there are other levels of granularity that could perform better (e.g., the optimized granularity from Nocera and Gardoni 2022). However, for wastewater networks, the results indicate that finer network granularity is more important than exactly representing the flow dynamics. Moving on to comparison between the full networks, model KF and model CF performed similarly. Because model KF took 46.7 seconds per simulation while model CF took 1.16 seconds per simulation, model KF is not preferred. The results indicate that model DF or model CF should be used. Model DF is preferred if predicting the exact value of the performance measure is needed. However, for many applications, the level of accuracy of model CF may be acceptable. Moreover, model CF gives the same spatial distribution of the performance measure as model DF. If priority is on locating vulnerable areas of a community, as opposed to estimating the performance measure exactly, model CF may be preferred. Moreover, given that a simulation of model CF takes approximately 0.6% the time of model DF, some optimization or uncertainty quantification tasks that are infeasible with model DF may be feasible with model CF. Moreover, the case study in this paper was for a small community. The running time of flow analysis using dynamic wave routing on full networks for large communities may be infeasibly large when many simulations are needed.

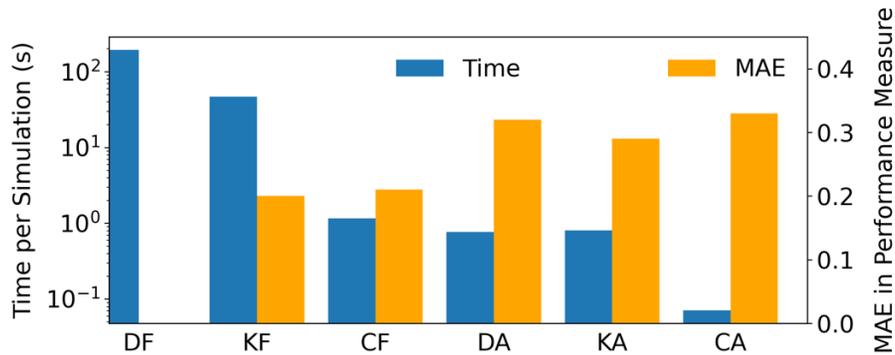

**Figure 9.** Accuracy and running time for the tested models

## 7. Conclusions

In this paper, we tested the performance of different representations of wastewater infrastructure for risk analysis. Specifically, we identified two levels of network granularity: the full network and the arterial network; and we identified three levels of fidelity of the governing equations: the dynamic wave equations, the kinematic wave equations, and connectivity analysis. Combinations of these two characteristics defined the 6 models that we tested. We simulated damage to pumps, tanks, and pipes based on assumed probabilities of damage states. We defined a performance measure for each model to represent the functionality of each node in the network given the damage. We demonstrated the model selection for the wastewater network of Seaside, Oregon. We compared the predictions from all of the models. As compared to the model with the full network that used the dynamic wave equations (i.e., the most accurate model), no other model had a mean absolute error in the nodal performance measure of less than 0.2. However, the spatial distribution of the performance measure was consistent across the models. The accuracy of the models that used the arterial network, even with the highest fidelity of the governing equations, was poorer than connectivity analysis with the full network. We conclude that, among the options considered, the full network should be used. We also conclude that the kinematic wave equations should not be used as they take significantly longer than connectivity analysis without providing any additional information. The kinematic wave equations do not capture surcharge or backwater effects in the flow, which are important physical phenomena for the performance of damaged wastewater collection systems. We conclude that the full network with either the dynamic wave equations or connectivity analysis should be used. If predicting the performance measure exactly is important, the dynamic wave equations are preferred. If the spatial distribution of damage is more important, and some error in the value of the performance measure is acceptable, connectivity analysis may be preferred.

**Repository**

The scripts used to generate results for this paper are available publicly at the following repository:

https://github.com/ajdunton/ww-model-selection